\documentclass[a4paper,reqno,12pt]{amsart}
\usepackage{verbatim}
\usepackage{amssymb,amsmath,amsthm}
\usepackage{amsfonts}
\usepackage{array}

\theoremstyle{remark}

\newcommand*\pFq[5]{{}_{#1}F_{#2}{\left[\genfrac..{0pt}{}{#3}{#4};#5\right]}}

\def\and{\quad\mbox{and}\quad}
\begin{document}

\setcounter{page}{1}

\title[]{A supercongruence involving\\
cubes of Catalan numbers}

\author{Roberto~Tauraso}

\address{Dipartimento di Matematica, 
Universit\`a di Roma ``Tor Vergata'', 
via della Ricerca Scientifica,
00133 Roma, Italy}
\email{tauraso@mat.uniroma2.it}




\begin{abstract} We mainly show a supercongruence for a truncated series with cubes of Catalan numbers which extends a result by Zhi-Wei Sun. 
\end{abstract}

\maketitle

\section{Introduction}
Let us consider the sum
$$\sum_{k=0}^n\left(\frac{C_k}{4^k}\right)^d$$
where $d$ is a positive integer and $C_k=\frac{1}{k+1}\binom{2k}{k}$ is the $k$-th Catalan number. For $d=1$ and $d=2$, it is easy to find a closed formula:
$$\sum_{k=0}^n\frac{C_k}{4^k}=2-\frac{\binom{2n+1}{n}}{4^n}\quad\mbox{and}\quad
\sum_{k=0}^n\frac{C_k^2}{16^k}=-4+(5+4n)\frac{\binom{2n+1}{n}^2}{16^n}.
$$
Thus, it follows that
$$\sum_{k=0}^{\infty}\frac{C_k}{4^k}=2\quad,\quad \sum_{k=0}^{\infty}\frac{C_k^2}{16^k}=-4+\frac{16}{\pi},$$
and, for any prime $p>2$,
$$\sum_{k=0}^{(p-1)/2}\frac{C_k}{4^k}\equiv_{p^2}2-2(-1)^{(p-1)/2}p\quad,\quad 
\sum_{k=0}^{(p-1)/2}\frac{C_k^2}{16^k}\equiv_{p^3}-4+12p^2$$
where we use the notation $a\equiv_m b$ to mean $a\equiv b \pmod{m}$.

When $d=3$ it seems that there is no closed formula for the partial sum. However, by Dixon's theorem \cite[p.13]{Ba35}, we are able to evaluate the infinite sum:
$$\sum_{k=0}^{\infty}\frac{C_k^3}{64^k}=
8\left(1 - \pFq{3}{2}{-\frac{1}{2},-\frac{1}{2},-\frac{1}{2}}{1,1}{1}\right)
=8-\frac{384\pi}{\Gamma^4(\frac{1}{4})}.$$
What about the related congruence?  We are going to show that for any prime $p>2$,
\begin{equation}\label{CCC3}
\sum_{k=0}^{(p-1)/2}\frac{C_k^3}{64^k}\equiv_{p^3} \left \{ \begin{array}{ll}\displaystyle 8-\frac{24p^2}{\Gamma_p^4\left(\frac 14\right)}   & \text{ if } p\equiv_4 1, \vspace{3mm}\\
 \displaystyle 8-\frac{384}{\Gamma_p^4\left(\frac 14\right)}  
& \text{ if } p\equiv_4 3. \end{array} \right.
\end{equation}
where $\Gamma_p$ is the Morita's $p$-adic Gamma function which is defined as the continuous extension to the set of all $p$-adic integers $\mathbb{Z}_p$ of the sequence
$$n\to (-1)^n\prod_{\substack{0\le k < n\\
(k,p)=1}}k$$
(see \cite[Chapter 7]{Ro00} for a detailed introduction to $\Gamma_p$). 
The above congruence modulo $p$ has been showed by Zhi-Wei Sun in \cite{Szw11}[Theorem 1.2]. Modulo $p^2$, the case $p\equiv_4 1$ is implied by  \cite{Szw12}[Theorem 1.3] (see also \cite{MP18}[Theorem 1.1]). 

\section{A bunch of identities}

The following one-parameter formula is the identity 6.34 in Gould's collection \cite{Go72} (for $j=0$ see \cite{Be30})
\begin{equation}\label{BELL}
\sum_{k=0}^nA(n,k,j)=\frac{1+(-1)^{n-j}}{2}\binom{n+j}{(n+j)/2}\binom{n-j}{(n-j)/2}
\end{equation}
where
$$A(n,k,j)=\binom{n}{k}\binom{n+k}{k}\binom{2k}{k+j}(-4)^{n-k}.$$
By using the partial fraction expansion we have that
$$\frac{A(n,k,1)}{A(n,k,0)}=1-\frac{1}{k+1},$$
and, by \eqref{BELL}, we get
\begin{equation}\label{ID1}
\sum_{k=0}^n\frac{\binom{n}{k}\binom{n+k}{k}\binom{2k}{k}}{(-4)^k(k+1)}=\frac{\binom{2\lfloor n/2\rfloor}{\lfloor n/2\rfloor}^2}{16^{\lfloor n/2\rfloor}}\cdot\begin{cases}
\displaystyle 1
&\text{if $n\equiv_2 0$},\vspace{3mm}\\
\displaystyle \frac{n}{n+1}
&\text{if $n\equiv_2 1$}.
\end{cases}
\end{equation}
In a similar way, the expansions
$$
\frac{A(n+1,k+1,1)-(n^2-n)A(n,k,2)}{A(n,k,0)}=
4+n-n^2+\frac{4n^2+4n-2}{k+1}-\frac{2n(n+1)}{(k+1)^2},$$
and
$$\frac{A(n+1,k+1,0)}{A(n,k,0)}=
4+\frac{8n+2}{k+1}+\frac{4n^2-2}{(k+1)^2}
-\frac{2n(n+1)}{(k+1)^3},
$$
yield respectively
\begin{equation}\label{ID2}
\sum_{k=0}^n\frac{\binom{n}{k}\binom{n+k}{k}\binom{2k}{k}}{(-4)^k(k+1)^2}=\frac{\binom{2\lfloor n/2\rfloor}{\lfloor n/2\rfloor}^2}{16^{\lfloor n/2\rfloor}}\cdot\begin{cases}
\displaystyle 2
&\text{if $n\equiv_2 0$},\vspace{3mm}\\
\displaystyle \frac{2n^2+2n-1}{(n+1)^2} 
&\text{if $n\equiv_2 1$},\end{cases}
\end{equation}
and
\begin{equation}\label{ID3}
\sum_{k=0}^n\frac{\binom{n}{k}\binom{n+k}{k}\binom{2k}{k}}{(-4)^k(k+1)^3}=-\frac{2}{n(n+1)}
+\frac{\binom{2\lfloor n/2\rfloor}{\lfloor n/2\rfloor}^2}{16^{\lfloor n/2\rfloor}}\cdot\begin{cases}
\displaystyle \frac{(2n+1)^2}{n(n+1)}
&\text{if $n\equiv_2 0$},\vspace{3mm}\\
\displaystyle \frac{4n^4+8n^3+3n^2-n+1}{n(n+1)^3}
&\text{if $n\equiv_2 1$}.
\end{cases}
\end{equation}
We would like to point out that, by the same approach, we  are able to find an explicit formula for
$$\sum_{k=0}^n\frac{Q(k)\binom{n}{k}\binom{n+k}{k}\binom{2k}{k}}{(-4)^k(k+1)^3}$$
where $Q$ is any polynomial in $\mathbb{Z}[x]$.
For example
\begin{equation}\label{ID4}
\sum_{k=0}^n\frac{k^3\binom{n}{k}\binom{n+k}{k}\binom{2k}{k}}{(-4)^k}=
\displaystyle\frac{\binom{2\lfloor n/2\rfloor}{\lfloor n/2\rfloor}^2}{16^{\lfloor n/2\rfloor}}\cdot\begin{cases}
\displaystyle \frac{n^2(n+1)^2(2n+1)^2}{15} 
&\text{if $n\equiv_2 0$},\vspace{3mm}\\
\displaystyle -\frac{n^2(4n^4+8n^3+3n^2-n+1)}{15}
&\text{if $n\equiv_2 1$}.
\end{cases}
\end{equation}
Moreover, by using \cite[Lemma 4.2]{Ta18}, we have that, for any positive even number $n$,
\begin{align}\label{ID5}
\sum_{k=0}^n\frac{\binom{n}{k}\binom{n+k}{k}\binom{2k}{k}H^{(2)}_{k+1}}{(-4)^k(k+1)^3}
&=\frac{\binom{n}{n/2}^2}{4^n}\left(16+\frac{(2n+1)^2\sum_{k=1}^n\frac{(-1)^k}{k^2}}{n(n+1)}\right)\nonumber\\
&\qquad\quad-\frac{4^n}{\binom{n}{n/2}^2}\cdot \frac{4n^4+8n^3+3n^2-n+1}{(n(n+1))^3}.
\end{align}
where $H_k^{(r)}=\sum_{j=1}^k1/j^r$ is the $k$-th harmonic number of order $r\geq 1$.

\section{Proof of \eqref{CCC3}}

First of all we need a stronger versions of \cite{Ta18}[Lemma 3]: for any prime $p>3$  
\begin{align}\label{G4plus}
\frac{\binom{2\lfloor n/2\rfloor}{\lfloor n/2\rfloor}^2}{16^{\lfloor n/2\rfloor}}
\equiv_{p^2}\begin{cases}
\displaystyle-\Gamma^4_p\left(\frac{1}{4}\right) 
\left(1+\frac{p^2E_{p-3}}{2}\right)
&\text{if $p\equiv_4 1$},\vspace{3mm}\\
\displaystyle\frac{16+32p+p^2(48-8E_{p-3})}{\Gamma^{4}_p\left(\frac{1}{4}\right)}
&\text{if $p\equiv_4 3$},
\end{cases}
\end{align}
where $n=(p-1)/2$ and $E_k$ is the $k$-th Euler number.
Indeed, for $p\equiv_4 3$,
\begin{align*}
\frac{\binom{2\lfloor n/2\rfloor}{\lfloor n/2\rfloor}^2}{16^{\lfloor n/2\rfloor}}
&=\frac{\Gamma_p^2\left(\frac{p+3}{4}\right)}{
\left(\frac{p-1}{4}\right)^2\Gamma_p^2\left(\frac{p+1}{4}\right)}
\equiv_{p^3}
\frac{\Gamma_p^2\left(\frac{3}{4}\right)\left(
1-\frac{(-1)^{n}}{16}p^2H^{(2)}_{\lfloor p/4\rfloor}\right)^2}{
\left(\frac{p-1}{4}\right)^2\Gamma_p^2\left(\frac{1}{4}\right)}
\nonumber \\
&\equiv_{p^3}
\displaystyle\frac{16+32p+p^2(48-8E_{p-3})}{\Gamma^{4}_p\left(\frac{1}{4}\right)}
\end{align*}
because, by \cite[(20)]{Le38}, $H^{(2)}_{\lfloor p/4\rfloor}\equiv_p (-1)^{n}4E_{p-3}$.
The proof for the case $p\equiv_4 1$ is similar:
\begin{align*}\frac{\binom{2\lfloor n/2\rfloor}{\lfloor n/2\rfloor}^2}{16^{\lfloor n/2\rfloor}}
&=
-\frac{\Gamma_p^2\left(\frac{p+1}{4}\right)}{
\Gamma_p^2\left(\frac{p+3}{4}\right)}
\equiv_{p^3}
-\frac{\Gamma_p^2\left(\frac{1}{4}\right)\left(
1+\frac{(-1)^{n}}{16}p^2H^{(2)}_{\lfloor p/4\rfloor}\right)^2}{
\Gamma_p^2\left(\frac{3}{4}\right)}
\nonumber \\
&\equiv_{p^3}
\displaystyle-\Gamma^4_p\left(\frac{1}{4}\right) 
\left(1+\frac{p^2E_{p-3}}{2}\right).
\end{align*}
For $0\leq k\leq n$,
$$
\binom{n}{k}\binom{n+k}{k}(-1)^k=\binom{2k}{k}\frac{\prod_{j=1}^k((2j-1)^2-p^2)}{4^k(2k)!}
\equiv_{p^3}
\frac{\binom{2k}{k}^2}{16^k}\left(1-p^2\sum_{j=1}^k\frac{1}{(2j-1)^2}\right).
$$
Hence, by identity \eqref{ID3}, congruence \eqref{CCC3} is implied by 
\begin{equation}\label{CCCH}
\sum_{k=0}^{(p-1)/2}\frac{C_k^3}{64^k}\sum_{j=1}^k\frac{1}{(2j-1)^2}\equiv_{p} \left \{ \begin{array}{ll} 
 \displaystyle -8-\frac{24}{\Gamma_p^4\left(\frac 14\right)}-4\Gamma_p^4\left(\frac 14\right)  
& \text{ if } p\equiv_4 1, \vspace{3mm}\\
\displaystyle -8+\frac{192(5-E_{p-3})}{\Gamma_p^4\left(\frac 14\right)} & \text{ if } p\equiv_4 3
. \end{array} \right.
\end{equation}
In order to show \eqref{CCCH}, notice that
 $$\frac{C_k}{4^k}\equiv_p \frac{(-1)^{k}}{n+1}\binom{n+1}{k+1}\quad\text{and}\quad
 \sum_{j=1}^k\frac{1}{(2j-1)^2}\equiv_p -\frac{H_{n-k}^{(2)}}{4}.$$
For $p\equiv_4 3$ we may use the identity \cite{Sh07}{(1.12)},
$$\sum_{k=0}^{2m}(-1)^k\binom{2m}{k}^3H_k^{(2)}=\frac{(-1)^m}{2}\,\frac{(3m)!}{(m!)^3}\,(H_m^{(2)}+H_{2m}^{(2)})$$
with $2m=n+1$. Then
\begin{align*}
\sum_{k=0}^{(p-1)/2}\frac{C_k^3}{64^k}\sum_{j=1}^k\frac{1}{(2j-1)^2}
&\equiv_p \frac{1}{4(n+1)^3}\sum_{k=0}^{n}(-1)^{k+1}\binom{n+1}{k+1}^3H_{n+1-(k+1)}^{(2)}
\\
&=\frac{1}{4(n+1)^3}\left((-1)^{n+1}\sum_{k=0}^{n}(-1)^{k+1}\binom{n+1}{k+1}^3H_{k+1}^{(2)}
-H_{n+1}^{(2)}\right)\\
&\equiv_p 2\sum_{k=1}^{2m}(-1)^{k}\binom{2m}{k}^3H_{k}^{(2)}-2H_{2m}^{(2)}
=(-1)^m\,\frac{(3m)!}{(m!)^3}\,(H_m^{(2)}+H_{2m}^{(2)})-2H_{2m}^{(2)}\\
&\equiv_p \frac{192(5-E_{p-3})}{\Gamma_p^4\left(\frac 14\right)}-8
\end{align*}
where $(-1)^m\,(3m)!/(m!)^3\equiv_p 48/\Gamma_p^4(1/4)$,
$H_{m}^{(2)}=H^{(2)}_{\lfloor p/4\rfloor}+1/(\frac{p+1}{4})^2\equiv_p (-1)^{n}4E_{p-3}+16$ and $H_{2m}^{(2)}=H_{n}^{(2)}+1/(\frac{p+1}{2})^2\equiv_p 0+4$.

\medskip

\noindent We have to deal with the case $p\equiv_4 1$ differently, because, for $2m=n$, the sum
$$\sum_{k=0}^{2m+1}(-1)^k\binom{2m+1}{k}^3H_k^{(2)}$$
cannot be expressed as a closed formula involving harmonic numbers (see \cite{Sh07}{(1.15)}). 
By \eqref{ID5} and \eqref{G4plus},
\begin{align*}
\sum_{k=0}^{(p-1)/2}\frac{C_k^3}{64^k}\sum_{j=1}^k\frac{1}{(2j-1)^2}
&\equiv_p \frac{1}{4(n+1)^3}\sum_{k=0}^{n}(-1)^{k+1}\binom{n+1}{k+1}^3H_{n+1-(k+1)}^{(2)}
\\
&=\frac{(-1)^{n}}{4}\sum_{k=0}^{n}(-1)^{k}\binom{n}{k}^3\frac{H_{k+1}^{(2)}}{(k+1)^3}
-\frac{H_{n+1}^{(2)}}{4(n+1)^3}\\
&\equiv_p \frac{1}{4}\sum_{k=0}^n\frac{\binom{n}{k}\binom{n+k}{k}\binom{2k}{k}H^{(2)}_{k+1}}{(-4)^k(k+1)^3}-2H_{n+1}^{(2)}\\
&=\frac{\binom{n}{n/2}^2}{4^n}\left(4+\frac{(2n+1)^2\sum_{k=1}^n\frac{(-1)^k}{k^2}}{4n(n+1)}\right)\nonumber\\
&\qquad\quad-\frac{4^n}{\binom{n}{n/2}^2}\cdot \frac{4n^4+8n^3+3n^2-n+1}{4(n(n+1))^3}
-2H_{n+1}^{(2)}\\
&\equiv_p -4\Gamma_p^4\left(\frac 14\right)-\frac{24}{\Gamma_p^4\left(\frac 14\right)}-8 
\end{align*}
and the proof is complete.

\medskip

As a final remark, we notice that, by using a similar approach, starting from identity \eqref{ID4}, it follows that for any prime  $p>5$,
\begin{equation}\label{CCC2}
\sum_{k=0}^{(p-1)/2}\frac{k^3\binom{2k}{k}^3}{64^k}\equiv_{p^3} \left \{ \begin{array}{ll} 
 \displaystyle -\frac{p^2}{40\Gamma_p^4\left(\frac 14\right)} & \text{ if } p\equiv_4 1, \vspace{3mm}\\
 \displaystyle -\frac{2}{5\Gamma_p^4\left(\frac 14\right)}  
& \text{ if } p\equiv_4 3. \end{array} \right.
\end{equation}
This congruence modulo $p$ appeared in \cite{Szw11}[Theorem 1.2] 
whereas,  the case $p\equiv_4 1$ modulo $p^2$ is implied by \cite{Szw12}[Theorem 1.3].

\end{document}